\documentclass[10pt]{article}
\usepackage{amsmath,amsfonts}
\normalbaselines
\newtheorem{teor}{Theorem}
\newtheorem{prop}[teor]{Proposition}
\newtheorem{lema}[teor]{Lemma}
\newtheorem{coro}[teor]{Corollary}
\newtheorem{rem}[teor]{Remark}

\newtheorem{ejem}[teor]{Example}

\newenvironment{demo}{\rm \trivlist \item[\hskip \labelsep{\it
      Proof}.]}{\nopagebreak \hfill $\square$ \endtrivlist}

\title{Uniqueness of complete maximal hypersurfaces in spatially open $(n+1)$-dimensional Robertson-Walker spacetimes\\ with flat fiber}

\author{Jos\'e A. S. Pelegr\'in${}^{a}$, Alfonso
Romero${}^{a}$ and Rafael M. Rubio${}^{b}$ \\[6mm]
${}^a$ Departamento de Geometr\'\i a y Topolog\'\i a, \\ [0.5mm]
Universidad de Granada, 18071 Granada, Spain \\ E-mails\textup{:
\texttt{jpelegrin@ugr.es, aromero@ugr.es}} \\[3mm]
${}^b$ Departamento de Matem\'aticas, Campus de Rabanales, \\[0.5mm] Universidad de
C\'ordoba, 14071 C\'ordoba, Spain,\\[0.5mm] E-mail\textup{: \texttt{rmrubio@uco.es}}\\[3mm]}

\date{}

\begin{document}

\maketitle

\thispagestyle{empty}

\begin{center}
Dedicated to Professor \'Angel Ferr\'andez
\end{center}

\begin{abstract}
In this paper,  under natural geometric and physical
assumptions we provide new uniqueness and non-existence
results for complete maximal hypersurfaces in spatially open Robertson-Walker spacetimes whose fiber is flat. Moreover, our results are applied to relevant spacetimes as the steady state spacetime, Einstein-de Sitter spacetime and radiation models.
\end{abstract}
\vspace*{5mm}




\section{Introduction}
The importance in General Relativity of maximal and constant mean
curvature spacelike hypersurfaces in spacetimes is well-known; a
summary of several reasons justifying it can be found in \cite{M-T}. Each maximal hypersurface can describe, in some relevant cases, the transition between the expanding and 
contracting phases of a relativistic universe. Moreover, the existence of constant mean curvature (and in particular maximal) hypersurfaces 
is necessary for the study of the structure of singularities in the space of solutions
to the Einstein equations. Also, the deep understanding of this kind of hypersurfaces is essential to prove the positivity of the gravitational mass.
They are also interesting for Numerical Relativity, where maximal hypersurfaces are used for integrating forward in time. All these physical aspects can be found in 
\cite{JKG}, \cite{M-T} and references therein.

A maximal hypersurface is (locally) a critical point for 
a natural variational problem, namely of the area functional (see, for
instance, \cite{brasil}). From a mathematical point of view, it is necessary to study the maximal hypersurfaces of a spacetime in order to understand its structure \cite{Bar}. Especially,
for some asymptotically flat spacetimes, the existence of a foliation by maximal hypersurfaces is established (see, for instance,
\cite{BF} and references therein). The existence results and, consequently, uniqueness appear as kernel topics.
Classical papers dealing with uniqueness results for constant mean curvature (CMC) hypersurfaces are
\cite{BF}, \cite{Ch} and \cite{M-T}, although a previous relevant
result in this direction was the proof of the Bernstein-Calabi
conjecture \cite{Ca} for the $n$-dimensional Lorentz-Minkowski
spacetime given by Cheng and Yau \cite{CY}. In \cite{BF},
Brill and Flaherty replaced the Lorentz-Minkowski spacetime by a
spatially closed universe, and proved uniqueness results for CMC hypersurfaces in the large by
assuming $\overline{{\mathrm{Ric}}}(z,z)>0$ for every timelike vector $z$. In
\cite{M-T}, this energy condition was relaxed by Marsden and Tipler
to include, for instance, non-flat vacuum spacetimes. More recently, 
Bartnik proved in \cite{Bar} very general existence theorems and
consequently averred that it would be useful to find new
satisfactory uniqueness results. Still more recently, in
\cite{A-R-S1} Al\'ias, Romero and S\'anchez proved new uniqueness
results in the class of spacetimes that they called spatially closed Generalized Robertson-Walker (GRW) spacetimes
(which includes the spatially closed Robertson-Walker spacetimes)
under the Timelike Convergence Condition. This GRW spacetimes differ from the classical Robertson-Walker spacetimes due to the fact that, despite being both defined as the warped product of an open interval endowed with negative definite metric and a Riemannian manifold as a fiber, they do not necessarily have constant sectional curvature. Finally, 
Romero, Rubio and Salamanca provided uniqueness results, in the
maximal case, for spatially parabolic GRW
spacetimes in \cite{RRS}, which are spatially open models whose fiber is a parabolic
Riemannian manifold. Moreover, making use of a well-known generalized maximum principle, the same authors obtain in \cite{RRS2} new uniqueness results in  other relevant spatially open GRW spacetimes for complete maximal hypersurfaces which are between two spacelike slices (time bounded) and have a bounded hyperbolic angle.

 In this paper we focus on the problems of uniqueness and non-existence of complete maximal hypersurfaces immersed in a spatially open Robertson-Walker spacetime with flat fiber. Note that these models have aroused a great deal of interest, since recent observations have shown that the current universe is very close to a spatially flat geometry \cite{CST}. This is actually a natural result from inflation in the early universe \cite{Lid}. We will give results that can be used when the fiber is $\mathbb{R}^n$, which is not parabolic for $n \geq 3$ and therefore, cannot be studied in  arbitrary dimension using previous methods. What is more important, we will not need the hyperbolic angle of the hypersurface to be bounded, which was a restrictive assumption used in previous works studying the spatially open case. Since we are not imposing this restriction, we are able to deal with spacelike hypersurfaces approaching the null Scri boundary at infinity, such as hyperboloids in Minkowski spacetime.

Our paper is organized as follows. Section \ref{s2} is devoted to introduce the basic notation used to describe spacelike hypersurfaces in GRW spacetimes. In Section \ref{tr} we provide an inequality involving the hyperbolic angle of a maximal hypersurface immersed in a GRW spacetime whose fiber is Ricci-flat and obeys the Null Convergence Condition (see Lemma \ref{lemachulo}). This inequality will play a crucial role in our results. In Section \ref{mr} we obtain a uniqueness result for complete maximal hypersurfaces (Theorem \ref{teomaxi}). In order to obtain it, the fundamental tool will be a Liouville-type theorem applied to the inequality obtained in  Lemma \ref{lemachulo}. Finally, in Subsection \ref{phy} we give several interpretations of our mathematical assumptions, their physical meaning and their compatility with the usual Energy Conditions using a perfect fluid model.

\section{Preliminaries}
\label{s2} \noindent
Let $(F,g_{_F})$ be an $n(\geq 2)$-dimensional (connected)
Riemannian manifold, $I$ an open interval in $\mathbb{R}$ endowed with the metric $-dt^2$  and $f$ a
positive smooth function defined on $I$. Then, the product manifold
$I \times F$ endowed with the Lorentzian metric
\begin{equation}\label{metrica}
\bar{g} = -\pi^*_{_I} (dt^2) +f(\pi_{_I})^2 \, \pi_{_F}^* (g_{_F})
\, ,
\end{equation}
where $\pi_{_I}$ and $\pi_{_F}$ denote the projections onto $I$ and
$F$, respectively, is called a \emph{GRW spacetime} with
\emph{fiber} $(F,g_{_F})$, \emph{base} $(I,-dt^2)$ and \emph{warping
function} $f$. If the fiber has constant sectional curvature, it is called a \emph{Robertson-Walker spacetime}.

In any GRW spacetime $\overline{M}=I\times_f F$, the
coordinate vector field $\partial_t:=\partial/\partial t$ is
(unitary) timelike, and hence $\overline{M}$ is time-orientable. The vector field $\partial_t$ plays a key role in the study of these spacetimes, since it constitutes a proper time synchronizable reference frame, which is geodesic, spatially conformal and irrotational \cite{S}. On the other hand, if we consider the timelike vector field $K: =~ f({\pi}_I)\,\partial_t$, from the relation between the
Levi-Civita connection of $\overline{M}$ and those of the base and
the fiber \cite[Cor. 7.35]{O'N}, it follows that
\begin{equation}\label{conexion} \overline{\nabla}_XK =
f'({\pi}_I)\,X
\end{equation}
for any $X\in \mathfrak{X}(\overline{M})$, where $\overline{\nabla}$
is the Levi-Civita connection of the Lorentzian metric
(\ref{metrica}). Thus, $K$ is conformal and its metrically equivalent $1$-form is closed.

From (\ref{conexion}) we easily see that the divergence on $\overline{M}$ of the reference frame $\partial_t$ satisfies 
${\rm div}(\partial_t) = n \frac{f'(t)}{f(t)}$. Therefore, the observers in $\partial_t$ are spreading out (resp. coming 
together) if $f'>0$ (resp. $f'<0$).

Given an $n$-dimensional manifold $M$, an immersion $\psi: M
\rightarrow \overline{M}$ is said to be \emph{spacelike} if the
Lorentzian metric (\ref{metrica}) induces, via $\psi$, a Riemannian
metric $g_{_M}$ on $M$. In this case, $M$ is called a spacelike
hypersurface. We will denote by $\tau:=\pi_I\circ \psi$ the
restriction of $\pi_I$ along $\psi$.

The time-orientation of $\overline{M}$ allows to take, for each
spacelike hypersurface $M$ in $\overline{M}$, a unique unitary
timelike vector field $N \in \mathfrak{X}^\bot(M)$ globally defined
on $M$ with the same time-orientation as $\partial_t$, i.e., such
that $\bar{g}(N,\partial_t)\leq -1$ and $\bar{g}(N,\partial_t)=
-1$ at a point $p\in M$ if and only if $N = \partial_t$ at
$p$. We will denote by $A$ the shape operator associated to $N$.
Then, the \emph{mean curvature function} associated to $N$ is given
by $H:= -(1/n) \mathrm{trace}(A)$. As it is well-known, the mean
curvature is constant if and only if the spacelike hypersurface is,
locally, a critical point of the $n$-dimensional area functional for
compactly supported normal variations, under certain constraints of
the volume. When the mean curvature vanishes identically, the
spacelike hypersurface is called a \emph{maximal} hypersurface.

For a spacelike hypersurface $\psi: M \rightarrow \overline{M}$ with
Gauss map $N$, the \emph{hyperbolic angle} $\varphi$, at any point
of $M$, between the unit timelike vectors $N$ and $\partial_t$, is
given by $\cosh \varphi=-\bar{g}(N,\partial_t)$. For simplicity,
throughout this paper we will refer to $\varphi$  as the
\emph{hyperbolic angle function} on $M$.

In any GRW spacetime $\overline{M}= I \times_f F$ there
is a remarkable family of spacelike hypersurfaces, namely its
spacelike slices $\{t_{0}\}\times F$, $t_{0}\in I$. It can be easily
seen that a spacelike hypersurface in $\overline{M}$ is a (piece of)
spacelike slice if and only if the function $\tau$ is constant.
Furthermore, a spacelike hypersurface in $\overline{M}$ is a (piece
of) spacelike slice if and only if the hyperbolic angle $\varphi$
vanishes identically. The shape operator of the spacelike slice
$\tau=t_{0}$ is given by $A=-f'(t_{0})/f(t_{0})\mathbb{I}$, where $\mathbb{I}$
denotes the identity transformation, and therefore its (constant)
mean curvature is $H=f'(t_{0})/f(t_{0})$. Thus, a spacelike slice
is maximal if and only if $f'(t_{0})=0$ (and hence, totally
geodesic).

\section{Set up} \label{tr}
Let $\psi: M \rightarrow \overline{M}$ be an $n$-dimensional
spacelike hypersurface immersed in a GRW spacetime
$\overline{M}= I \times_f F$. If we denote by
\[
\partial_t^T:= \partial_t+\overline{g}(N,\partial_t)N
\]
the tangential component of $\partial_t$ along $\psi$, then it is
easy to check that the gradient of $\tau$ on $M$ is
\begin{equation}\label{part}
\nabla \tau=-\partial_t^T
\end{equation}
and so
\begin{equation}\label{sinh}
|\nabla \tau|^2=g_{_M}(\nabla \tau,\nabla \tau)=\sinh^2 \varphi.
\end{equation}

Moreover, since the tangential
component of $K$ along $\psi$ is given by $K^T=K+\overline{g}(K,N)N$, a direct computation from
(\ref{conexion}) gives
\begin{equation}\label{gradcosh}
\nabla \overline{g}(K,N)=-AK^T
\end{equation}
where we have used (\ref{part}), and also
\[
\nabla \cosh \varphi=A\partial_t^T-\frac{f'(\tau)}{f(\tau)}\overline{g}(N,\partial_t)\partial_t^T.
\]

On the other hand, if we represent by $\nabla$ the Levi-Civita
connection of the metric $g_{_M}$, then the Gauss and Weingarten
formulas for the immersion $\psi$ are given, respectively, by
\begin{equation}\label{GF}
\overline{\nabla}_X Y=\nabla_X Y-g_{_M}(AX,Y)N
\end{equation}
and
\begin{equation}\label{WF}
AX=-\overline{\nabla}_X N,
\end{equation}
where $X,Y\in\mathfrak{X}({M})$. Then,  taking the tangential component in
(\ref{conexion}) and using (\ref{GF}) and (\ref{WF}), we get
\begin{equation}\label{KT}
\nabla_X K^T=-f(\tau)\overline{g}(N,\partial_t)AX+f'(\tau)X
\end{equation}
where $X\in\mathfrak{X}({M})$ and $f'(\tau):=f'\circ \tau$. Since
also $K^T=f(\tau)\partial_t^T$, it follows from (\ref{part}) and
(\ref{KT}) that the Laplacian of $\tau$ on $M$ is
\begin{equation}\label{laptau}
\Delta \tau=-\frac{f'(\tau)}{f(\tau)}\{n+|\nabla
\tau|^2\}-nH\overline{g}(N,\partial_t).
\end{equation}

Consequently,
\begin{eqnarray}
\Delta f(\tau) & = & f'(\tau) \Delta \tau+f''(\tau)|\nabla \tau|^2
\nonumber \\
& = & -\frac{f'(\tau)^2}{f(\tau)} \, n + |\nabla \tau|^2 f(\tau)
(\log f)''(\tau)+n H f'(\tau)\cosh \varphi \label{lapftau}
\end{eqnarray}
and so
\begin{eqnarray}
\Delta \left( f(\tau) \cosh \varphi\right)
& = &
\cosh \varphi \, \, \Delta f(\tau)+f(\tau) \, \, \Delta \cosh \varphi+2g_{_M}(\nabla f(\tau), \nabla \cosh \varphi) \nonumber \\
& = & -\frac{f'(\tau)^2}{f(\tau)} \, n\cosh \varphi +  f(\tau) \cosh \varphi \, \, \sinh^2 \varphi \, (\log f)''(\tau)+n H f'(\tau)\cosh^2 \varphi \nonumber \\
& & +f(\tau) \, \, \Delta \cosh \varphi-2 f'(\tau) g_{_M}(A\partial_t^T,
\partial_t^T)-2\frac{f'(\tau)^2}{f(\tau)}\, \, \cosh \varphi \, \,
\sinh^2 \varphi, \label{lapfcosh}
\end{eqnarray}
where we have used (\ref{part})-(\ref{gradcosh}).

Furthermore, if we assume that $M$ is a maximal hypersurface, we
get from the Codazzi equation for $M$ that (see \cite[Eq.
8]{A-R-S1})
\begin{equation}\label{ARS1995}
\Delta \left( f(\tau) \cosh \varphi\right)=-\Delta
\overline{g}(K,N)=-\overline{{\rm
Ric}}(K^T,N)+f(\tau)\cosh\varphi \,\, {\rm trace}(A^2)
\end{equation}
where $\overline{{\rm Ric}}$ stands for the Ricci tensor on
$\overline{M}$. Therefore, from (\ref{lapfcosh}) and (\ref{ARS1995})
we have
\begin{eqnarray}
\overline{{\rm Ric}}(K^T,N) & = &  f(\tau)\cosh\varphi \,\, {\rm trace}(A^2) +\frac{f'(\tau)^2}{f(\tau)} \, n\cosh \varphi \nonumber \\
& &  -  f(\tau) \cosh \varphi \, \, \sinh^2 \varphi \, (\log f)''(\tau)-f(\tau) \, \, \Delta \cosh \varphi \nonumber \\
& & +2 f'(\tau) g_{_M}(A\partial_t^T,
\partial_t^T)+2\frac{f'(\tau)^2}{f(\tau)}\, \, \cosh \varphi \, \,
\sinh^2 \varphi. \label{Ric1}
\end{eqnarray}

If we put $N=N_{_F}-\overline{g}(N,\partial_t)\partial_t$, where
$N_{_F}$ denotes the projection of $N$ on the fiber $F$, it is easy
to obtain from (\ref{metrica}) that
\begin{equation}\label{senh2}
\sinh^2\varphi=f(\tau)^2 \,\, g_{_F}(N_{_F},N_{_F}).
\end{equation}

Besides, from \cite[Cor. 7.43]{O'N} we know that
\begin{equation}\label{ricpartial}
\overline{{\rm
Ric}}(\partial_t,\partial_t)=-n\frac{f''(\tau)}{f(\tau)}
\end{equation}
and 
\begin{equation}\label{ricNF}
\overline{{\rm Ric}}(N_{_F},N_{_F})=\sinh^2\varphi
\left(\frac{f''(\tau)}{f(\tau)}+(n-1)\frac{f'(\tau)^2}{f(\tau)^2}
\right)
\end{equation}
where we have
used (\ref{senh2}) and the fact that $F$ is Ricci-flat. Then, from (\ref{ricpartial}) and (\ref{ricNF})
we obtain
\begin{eqnarray}
\overline{{\rm Ric}}(K^T,N) & = & -f(\tau)\cosh\varphi \,\, \overline{{\rm Ric}}(N_{_F},N_{_F})-f(\tau)\cosh \varphi \,\, \sinh^2 \varphi \,\, \overline{{\rm Ric}}(\partial_t,\partial_t) \nonumber \\
& = & (n-1) f(\tau)\cosh\varphi\,\,\sinh^2\varphi \,\,(\log f)''(\tau).
\label{Ric2}
\end{eqnarray}

Finally, from (\ref{Ric1}) and (\ref{Ric2}) we get
\begin{eqnarray}
\Delta \cosh \varphi & = & -n \cosh\varphi \sinh^2\varphi \,\,(\log f)''(\tau) + \frac{f'(\tau)^2}{f(\tau)^2}\, \, \cosh \varphi \left(n+2\sinh^2
\varphi\right) \nonumber \\
&& +\cosh\varphi \,\, {\rm
trace}(A^2)+2\frac{f'(\tau)}{f(\tau)}\,\,g_{_M}(A\partial_t^T,
\partial_t^T). \label{wegotit}
\end{eqnarray}

On the other hand, the square algebraic trace-norm of the Hessian
tensor of $\tau$ is just
\[
\vert {\rm Hess}(\tau)\vert^2={\rm trace} (H_\tau\circ H_\tau),
\]
where $H_\tau$ denotes the operator defined by
$g_{_M}(H_\tau(X),Y):={\rm Hess}(\tau)(X,Y)$ for all $X,Y\in\mathfrak{X}({M})$.

Taking the tangential component in (\ref{conexion}) and using
(\ref{part}) we get that
\begin{eqnarray}
\vert{\rm Hess}(\tau)\vert^ 2 & = & \frac{f'(\tau)^ 2}{f(\tau)^
2}\left(n-1+\cosh^ 4\varphi\right)+\cosh^ 2\varphi\,{\rm
trace}(A^ 2) \nonumber \\
&& +2\frac{f'(\tau)}{f(\tau)}\cosh\varphi \, g_{_M}(A\partial_t^
T,\partial_t^ T). \label{hess}
\end{eqnarray}
Since $\vert{\rm Hess}(\tau)\vert^ 2\geq 0$, it is a straightforward
computation to obtain, making use of (\ref{wegotit}) and
(\ref{hess}), that
\begin{eqnarray}
\cosh \varphi\,\, \Delta \cosh \varphi & \geq &
- n \cosh^2\varphi \sinh^2\varphi \,\,(\log f)''(\tau) +n\frac{f'(\tau)^2}{f(\tau)^2}\cosh^2\varphi
\nonumber \\
&&
\hspace{-0.4cm}+2\frac{f'(\tau)^2}{f(\tau)^2}\cosh^2\varphi\sinh^2\varphi-
\frac{f'(\tau)^2}{f(\tau)^2}\left(n-1+\cosh^4\varphi\right).
\label{wegotit2}
\end{eqnarray}

Now, from (\ref{wegotit2}) we have

\begin{eqnarray}
\cosh \varphi\,\, \Delta \cosh \varphi & \geq & 
- n \frac{f''(\tau)}{f(\tau)} \cosh^2\varphi\sinh^2\varphi + n \frac{f'(\tau)^2}{f(\tau)^2} \cosh^2\varphi\sinh^2\varphi + n \frac{f'(\tau)^2}{f(\tau)^2} \cosh^2\varphi 
\nonumber \\
&& +2 \frac{f'(\tau)^2}{f(\tau)^2} \cosh^2\varphi\sinh^2\varphi - (n-1) \frac{f'(\tau)^2}{f(\tau)^2} - \frac{f'(\tau)^2}{f(\tau)^2} \left( \sinh^2\varphi + 1 \right) \cosh^2\varphi  \nonumber \\
&&= - n \frac{f''(\tau)}{f(\tau)} \cosh^2\varphi\sinh^2\varphi + (n-1)\frac{f'(\tau)^2}{f(\tau)^2} \cosh^2\varphi\sinh^2\varphi  \nonumber \\
&&+ 2 \frac{f'(\tau)^2}{f(\tau)^2}\cosh^2\varphi\sinh^2\varphi + (n-1)\frac{f'(\tau)^2}{f(\tau)^2} \cosh^2\varphi - (n-1) \frac{f'(\tau)^2}{f(\tau)^2} .
\label{wegotit3}
\end{eqnarray}

Since $\cosh^2\varphi \geq 1$, we obtain from (\ref{wegotit3})

\begin{equation}
\label{wegotit4}
\cosh \varphi\,\, \Delta \cosh \varphi  \geq  - n \frac{f''(\tau)}{f(\tau)} \cosh^2\varphi\sinh^2\varphi + (n+1)\frac{f'(\tau)^2}{f(\tau)^2} \cosh^2\varphi\sinh^2\varphi .
\end{equation}

Now, we introduce an extra assumption with physical meaning. A spacetime $\overline{M}$ obeys the
\emph{Null Convergence Condition} if its Ricci tensor
$\overline{\rm Ric}$ satisfies $\overline{\rm Ric}(z,z) \geq 0$, for
all null vectors $z$. When $\overline{M}$ is a GRW spacetime $I \times_f F$ with Ricci-flat fiber, it is not difficult to see that this energy condition is satisfied if and only if $(\log f)''(t)\leq 0$. Using this assumption, in (\ref{wegotit4}) we have

\begin{equation}
\label{wegotit5}
\cosh \varphi\,\, \Delta \cosh \varphi  \geq \left( (n+1)\frac{f'(\tau)^2}{f(\tau)^2} - n \frac{f''(\tau)}{f(\tau)} \right) \sinh^4 \varphi .
\end{equation}

Moreover, we know that

\begin{equation}
\label{lapsinh}
\frac{1}{2} \Delta \sinh^2 \varphi = \cosh \varphi \Delta \cosh \varphi + \vert \nabla \cosh \varphi \vert^2 \geq \cosh \varphi \Delta \cosh \varphi .
\end{equation}

From (\ref{wegotit5}) and (\ref{lapsinh}) we obtain the following result

\hyphenation{Ge-ne-ra-li-zed}

\begin{lema}
\label{lemachulo}
Let $\psi: M \rightarrow \overline{M}$ be an $n$-dimensional maximal hypersurface immersed in a GRW spacetime
$\overline{M}= I \times_f F$ with Ricci-flat fiber that obeys the Null Convergence Condition, then

\begin{equation}
\label{greatresult}
 \frac{1}{2} \Delta \sinh^2 \varphi \geq \left( (n+1)\frac{f'(\tau)^2}{f(\tau)^2} - n \frac{f''(\tau)}{f(\tau)} \right) \sinh^4 \varphi .
\end{equation}

\end{lema}

This inequality suggests us the use of the following lemma given by Nishikawa in \cite{N}, which extends and clarifies a technical step in Cheng and Yau's seminal paper \cite{CY}. In that paper, they used the lemma to study ${\rm trace}(A^2)$, whereas we will use it on the function $ \sinh^2 \varphi$.

\begin{lema}
\label{lema}{\rm \cite{N}}
Let $M$ be a complete Riemannian manifold whose Ricci curvature is bounded from below and let $u:M\longrightarrow\mathbb{R}$ be a non-negative smooth function on $M$. If there exists a constant $c>0$ such that $\Delta u\geq cu^2$, then $u$ vanishes identically on $M$.
\end{lema}

\begin{lema} \label{ricciacotado1}
Let $\psi: M \rightarrow \overline{M}$ be an $n$-dimensional maximal hypersurface immersed in a Robertson-Walker spacetime
$\overline{M}= I \times_f F$ with flat fiber that obeys the Null Convergence Condition.
Then, the Ricci 
curvature of $M$ must be non-negative.
\end{lema}

\begin{demo}
Given $p\in M$, let us take a local orthonormal frame $\left\{U_1,\ldots,U_n \right\}$ around $p$. From the
Gauss equation we get that the Ricci curvature of $M$, $\mathrm{Ric}$, satisfies

$$\mathrm{Ric}(Y,Y) \geq \sum_k \overline{g}(\overline{\mathrm{R}}(Y,U_k)U_k,Y),
$$ for all $Y \in \mathfrak{X}(M)$, where $\overline{\mathrm{R}}$ denotes the curvature tensor of $\overline{M}$ given by 

$$\overline{\mathrm{R}}(X,Y)Z =  \overline{\nabla}_X \overline{\nabla}_Y Z -\overline{\nabla}_Y \overline{\nabla}_X Z - \overline{\nabla}_{[X,Y]}Z, \quad \text{for all}\  X, Y, Z \in \mathfrak{X}(\overline{M}).$$

Now, from \cite[Prop. 7.42]{O'N} and using the fact that $F$ is flat, we have
\begin{eqnarray}
\label{ricbound}
\sum_k \overline{g}(\overline{\mathrm{R}}(Y,U_k)U_k,Y) &=&  (n-1) \frac{f'(\tau)^2}{f(\tau)^2} |Y|^2 -(n-2)(\log f)''(\tau) \, g(Y,\nabla \tau)^2  \nonumber \\
& & - (\log f)''(\tau) |\nabla \tau|^2 |Y|^2 .
\end{eqnarray}

 From these equations,
taking into account the assumptions, we have the Ricci curvature of $M$ to be non-negative.
\end{demo}

\begin{prop}
\label{prodi}
There is no complete maximal hypersurface $M$ in a spatially open Robertson-Walker spacetime $\overline{M}= I \times_f F$ with flat fiber that obeys the Null convergence Condition such that the restriction of the expanding/contracting function ${\rm div}(\partial_t)$ to $M$ satisfies 
$\inf_M |{\rm div}(\partial_t)| > 0.$
\end{prop}

\begin{demo}
From Lemma \ref{ricciacotado1} we get for any maximal hypersurface $M$ in $\overline{M}$

\begin{equation}
\label{ricdiv}
{\rm Ric}(Y,Y) \geq \frac{n-1}{n^2} {\rm div}(\partial_t)\Big{|}_M^2 |Y|^2,
\end{equation}

\noindent for all $Y \in \mathfrak{X}(M)$. Now, using our assumptions and (\ref{ricdiv}) we obtain that the Ricci curvature of $M$ is bounded from below by a positive constant. If $M$ is complete the classical Bonnet-Myers Theorem  
ensures its compactness. However, this contradicts the fact that in a spatially open spacetime there are no compact spacelike hypersurfaces.
\end{demo}

Proposition \ref{prodi} enables us to obtain the following non-existence results in some well-known Robertson-Walker spacetimes.

\begin{coro}
\label{coroste}
There are no complete maximal hypersurfaces in the $(n+1)$-dimensional steady state spacetime $\mathbb{R} \times_{e^t} \mathbb{R}^n$.
\end{coro}

\begin{coro}
\label{coroein}
There are no complete maximal hypersurfaces  bounded away from future infinity in the $(n+1)$-dimensional Einstein-de Sitter spacetime $\mathbb{R}^+ \times_{t^{2/3}} \mathbb{R}^n$.
\end{coro}

This generalizes and improves Rubio's result in \cite{Rub} to the case of arbitrary dimension. Analogously, we have

\begin{coro}
\label{cororad}
There are no complete maximal hypersurfaces  bounded away from future infinity in the $(n+1)$-dimensional Roberson-Walker Radiation Model spacetime $\mathbb{R}^+ \times_{(2at)^ {1/2}} \mathbb{R}^n$, with $a > 0$.
\end{coro}

\section{Main result} \label{mr}

As a consequence of Lemma \ref{lemachulo} and the Liouville-type result given in Lemma \ref{lema}, we can state our principal result

\begin{teor}
\label{teomaxi}

Let $\overline{M}= I \times_f F$ be a Robertson-Walker spacetime with flat fiber that obeys the Null Convergence Condition.
Then, the only complete maximal hypersurfaces immersed in 
$\overline{M}$ satisfying $\inf \left\{ (n+1)\frac{f'(\tau)^2}{f(\tau)^2} - n \frac{f''(\tau)}{f(\tau)} \right\}>0$ are the spacelike slices $\{t_0\}\times F$ with $f'(t_0)=0$. 
\end{teor}

\begin{demo}
From Lemmas \ref{lemachulo}, \ref{lema} and \ref{ricciacotado1}  we obtain that the hyperbolic angle must vanish identically on the maximal hypersurface. 
\end{demo}

\begin{rem} \normalfont Observe that the assumption on the function $(n+1)\frac{f'(\tau)^2}{f(\tau)^2} - n \frac{f''(\tau)}{f(\tau)}$ defined on the hypersurface is scarcely restrictive, even if combined with the NCC. In fact, if we consider its extension  $(n+1)\frac{f'(t)^2}{f(t)^2} - n \frac{f''(t)}{f(t)}$ defined on the spacetime, we have from the NCC that $(n+1)\frac{f'(t)^2}{f(t)^2} - n \frac{f''(t)}{f(t)}=\frac{f'(t)^2}{f(t)^2}-n(\log f)''(t)\geq 0$.

However, if we assume that the warping function is defined on the largest possible domain, i.e., it is inextendible, we can find two cases where the required inequality on the infimum does not hold:

\begin{enumerate}
\item When both $f'$ and $f''$ vanish simultaneously at some point in $I=]a, b[$. This obviously happens in the Lorentz-Minkowski spacetime, where an analogous uniqueness result does not hold. Note that $\mathbb{L}^{n+1}$ is a vacuum solution. What is more, if there is real presence of matter in the spacetime we can discard this case.

\item If $\lim\limits_{t \to b} \frac{f'(t)^2}{f(t)^2} = \lim\limits_{t \to b} (\log f)''(t) = 0$. This is the case in the Einstein-de Sitter spacetime. Even more, the inequality will not hold either in the less realistic case where $\lim\limits_{t \to a} \frac{f'(t)^2}{f(t)^2} = \lim\limits_{t \to a} (\log f)''(t) = 0$.
\end{enumerate}
\end{rem}

\begin{rem}\normalfont
Furthermore, this theorem improves some previous uniqueness results for complete maximal hypersurfaces (see \cite{RRS} and \cite{RRS2}, for instance) without making restrictive assumptions on the maximal hypersurface such as having a bounded hyperbolic angle or lying between two spacelike slices. Moreover, we would just need the fiber to be Ricci-flat instead of flat if we knew beforehand that the Ricci curvature of every complete maximal hypersurface in the spacetime is bounded from below. In this way, we could extend our results to GRW spacetimes with Ricci-flat fiber.
\end{rem}

We will give now two models where Theorem \ref{teomaxi} holds.

\begin{ejem}
\label{examp}
\normalfont
Let us consider the Robertson-Walker spacetime $\overline{M} = \mathbb{R} \times_f \mathbb{R}^n$ with warping function $f(t) = e^{-t^2}$.
This spacetime obeys NCC, since $ (\log f)''(t)= -2$.  Moreover, any maximal hypersurface immersed in $\overline{M}$ satisfies

$$\inf \left\{ (n+1)\frac{f'(\tau)^2}{f(\tau)^2} - n \frac{f''(\tau)}{f(\tau)} \right\} = \inf \{ 2n + 4\tau^2 \} > 0.$$
 
\noindent Therefore, the only complete maximal hypersurface in $\overline{M}$ is the spacelike slice $\{0\} \times \mathbb{R}^n$. 

This spacetime models a relativistic universe without singularities (in the sense of \cite[Def. 12.16]{O'N}) that goes from an expanding phase to a contracting one. The physical space in this transition of phase is represented by the spacelike slice $\{0\} \times \mathbb{R}^n$.

\end{ejem}

\begin{ejem}
\label{ext2}
\normalfont
We obtain another example of a Robertson-Walker spacetime satisfying the assumptions in Theorem 4 by considering $\overline{M} = I \times_f \mathbb{R}^n$. Where $I = ]-a, a[$ and the warping function is $f(t) = \sqrt{a^2 - t^2}$, being $a$ a positive constant. Let us remark that this spacetime behaves like the Robertson-Walker model proposed by Friedmann with constant sectional curvature of the fiber equal to one (see \cite[Chap. 12]{O'N}), since it has a big bang singularity at $t=-a$ as well as a big crunch at $t= a$ \cite[Def. 12.16]{O'N}.

For this spacetime, $ (\log f)''(t)= -\frac{a^2+t^2}{(a^2-t^2)^2} \leq 0$, so it satisfies NCC. Furthermore, for every maximal hypersurface in $\overline{M}$ 

$$\inf \left\{ (n+1)\frac{f'(\tau)^2}{f(\tau)^2} - n \frac{f''(\tau)}{f(\tau)} \right\} = \inf \left\{ \frac{n(a^2 + \tau^2) + \tau^2}{(a^2 - \tau^2)^2} \right\} > 0.$$

\noindent Hence, the only complete maximal hypersurface in this spacetime is the spacelike slice $\{0\} \times \mathbb{R}^n$, which represents the physical space in the transition from an expanding phase of the spacetime to a contracting one.

\end{ejem}

\begin{rem}
\normalfont
Note that Corollaries \ref{coroste}, \ref{coroein} and \ref{cororad} can also be easily deduced from Theorem \ref{teomaxi}.
\end{rem}

\subsection{Physical meaning of our mathematical assumptions}\label{phy}

There is currently great interest in the study of General Relativity in arbitrary dimensions
due to several reasons, such as the creation of unified theories
or the methodological considerations associated with the
possibility of understanding general features for the simpler
(2+l)-dimensional models (see \cite{Sha} and references therein).

We can build a simple cosmological model that will help us to better understand the hypotheses assumed in this article. In order to do so, we will start with our manifold $\overline{M}= I \times_f \mathbb{R}^n$, where the lines $ I \times {p}$ will be the worldlines of the galactic flow. Furthermore, if for each $p\in \mathbb{R}^n$ we parametrize $I\times \{ p \}$ by
$\gamma_{_p}(t)=(t,p)$ we can define our \emph{galaxies} $\gamma_{_p}$ as the integral curves of the velocity vector field $\partial_t$. In particular,
the function $t$ is the common proper time of all galaxies.
By taking $t$ as a constant, we get the hypersurface
\[
M(t)=\{t\} \times \mathbb{R}^n=\{(t,p): p\in \mathbb{R}^n\}.
\]

The distance between two \emph{galaxies} $\gamma_{_p}$ and $\gamma_{_q}$ in $M(t)$ is $f(t)d(p,q)$,
where $d$ is the Riemannian distance in the fiber $\mathbb{R}^n$. In particular,
when $f$ has positive derivative the spaces $M(t)$ are expanding.
Moreover, if $f''>0$, Robertson-Walker spacetimes model  universes in
accelerated expansion.

Some Robertson-Walker spacetimes
satisfying the \emph{Null Convergence Condition} can be suitable modified  models of gravity.
For instance, the steady state spacetime verifies  $f'(t)=e^t> 0$  and so the spaces $M(t)$ are expanding. In addition, this expansion is accelerated since $f''(t)=e^t$. Therefore, the steady state spacetime constitutes
an accelerated expanding spacetime.

On the other hand, astronomical evidence indicates that the universe can be modeled (in
smoothed, averaged form) as a spacetime containing a perfect fluid
whose \emph{molecules} are the galaxies. Classically, the dominant
contribution to the energy density of the galactic fluid is the mass
of the galaxies, with a smaller pressure due mostly to radiations.
Nevertheless, over the 90's, evidences for the most striking
result in modern cosmology have been steadily growing, namely the
existence of a cosmological constant which is driving the current
acceleration of the universe as first observed in \cite{SP},
\cite{R}. Different models for dark energy cosmology and their
equivalences can be seen in \cite{B}. Note that a positive vacuum
energy density resulting from a cosmological constant implies a
negative pressure and vice versa.

Thus, it is natural that several exact solutions to the Einstein field
equation
\begin{equation}\label{EE}
\overline{{\rm Ric}}-\frac{1}{2}\overline{S} \thinspace \overline{g \vphantom{S}}=8\pi T
\end{equation}
have been obtained by considering a continuous distribution of
matter as a perfect fluid.

Recall that a \emph{perfect fluid} (see, for example, \cite[Def.
12.4]{O'N}) on a spacetime $\overline{M}$ is a triple $(U,\rho,
\mathfrak{p})$ where
\begin{enumerate}
 \item $U$ is a timelike future-pointing unit vector field on $\overline{M}$ called the \emph{flow
 vector field}.
 \item $\rho, \mathfrak{p} \in C^\infty (\overline{M})$ are, respectively, the \emph{energy density} and the
 \emph{pressure} functions.
 \item The \emph{stress-energy momentum tensor} is
 $$
 T= (\rho + \mathfrak{p} ) \, U^{\overline{\flat}} \otimes U^{\overline{\flat}} + \mathfrak{p} \, \overline{g} ,
 $$ where $\overline{g}$ is the metric of the spacetime $\overline{M}$.
\end{enumerate}

For an instantaneous observer $v$, the quantity $T(v,v)$ is
interpreted as the energy density, i.e., the mass-energy per unit of
volume, measured by this observer. For normal matter, this
quantity must  be non-negative, i.e., the tensor $T$ must obey
the \emph{weak energy condition}. It is easy to see that an exact solution to (\ref{EE})
for a stress-energy tensor which obeys the Weak Energy Condition
must satisfy the null energy condition, that is, $\overline{{\rm
Ric}}(z,z)\geq 0$ for every null vector $z$. Nevertheless, perfect fluids  can also
be used to model another scenarios of universes at the dark energy
dominated stage (see \cite{CST}).

 In the case of a Robertson-Walker spacetime with flat fiber which is filled with a perfect fluid, the density and pressure functions are are given by

 $$8\pi\rho=\frac{n(n-1)}{2}\frac{f'^2}{f^2}$$ 
 and $$8\pi\mathfrak{p}=-(n-1)\frac{f''}{f}-\frac{(n-1)(n-2)}{2}\frac{f'^2}{f^2} ,$$ 
 
 \noindent being $f$ the warping function. For this family of classical models with positive density and non-negative pressure ($\rho>0$, $\mathfrak{p}\geq 0$), every maximal hypersurface included in a region of the spacetime in which the stress-energy momentum tensor is far from zero satisfies the condition $\inf \left\{ (n+1)\frac{f'(\tau)^2}{f(\tau)^2} - n \frac{f''(\tau)}{f(\tau)} \right\}>0$, which is equivalent in our model to $\inf \left\{ \frac{8 \pi}{n-1} \left( \frac{n^2 + 2}{n} \rho + n \mathfrak{p} \right) \right\} > 0 $. For instance, this happens near a physical singularity.

On the other hand, as the example of the steady state spacetime shows, the condition
 \par\noindent $\inf \left\{ (n+1)\frac{f'(\tau)^2}{f(\tau)^2} - n \frac{f''(\tau)}{f(\tau)} \right\}>0$ holds in certain models with negative pressure. 
 
Finally, we may wonder whether our assumptions in Theorem \ref{teomaxi} are compatible with other usual energy conditions \cite[Sec. 4.3]{HE}. In order to better understand these conditions, we will express them for a perfect fluid in terms of the energy density and pressure functions as well as using the warping function $f$ of the Robertson-Walker spacetime with flat fiber. Hence, for a Robertson-Walker spacetime with flat fiber filled with a perfect fluid we have

\begin{itemize}
\item The \emph{Weak Energy Condition} implies that $\rho \geq 0$ and $\rho + \mathfrak{p} \geq 0$; which is equivalent to \\ $(\log f)'' \leq 0$, since in our perfect fluid model the energy density is always non-negative.

\item The \emph{Strong Energy Condition} stipulates that $\rho + \mathfrak{p} \geq 0$ and $\rho + n \mathfrak{p} \geq 0$; or equivalently, $(\log f)'' \leq 0$ and $\frac{n-3}{2}\frac{f'^2}{f^2} + \frac{f''}{f} \leq 0$ (for $n=3$, this last equation leads to $f'' \leq 0$, which is equivalent to the \emph{Timelike Convergence Condition} in this case).

\item The \emph{Dominant Energy Condition} is written as $\rho \geq |\mathfrak{p}|$, which for positive pressure can be written as $(n-1)\frac{f'^2}{f^2} + \frac{f''}{f} \geq 0$ and for negative pressure implies $(\log f)'' \leq 0$.

\end{itemize}

Thus, our assumptions can be compatible in many cases with these energy conditions.

\section*{Acknowledgments} The authors would like to thank the referee for his deep reading and valuable suggestions. This paper has been partially supported by Spanish MINECO and ERDF project MTM2013-47828-C2-1-P.


\begin{thebibliography}{99}

\bibitem{A-R-S1} L.J. Al\'\i as, A. Romero and M. S\'anchez,
Uniqueness of complete spacelike hypersurfaces of constant mean
curvature in Generalized Robertson-Walker spacetimes, \emph{Gen. Relat. Gravit.}, \textbf{27} (1995), 71--84.


\bibitem{B} K. Bamba, S. Capozziello, S. Nojiri and S. D. Odintsov,
Dark energy cosmology: the equivalent description via different theoretical models and cosmography tests,
\emph{Astrophys. Space Sci.}, \textbf{342} (2012), 155--314.

\bibitem{Bar} R. Bartnik,
Existence of maximal surfaces in asymptotically flat spacetimes,
\emph{Commun. Math. Phys.}, \textbf{94} (1984), 155--175.



\bibitem{BF} D. Brill and F. Flaherty,
Isolated maximal surfaces in spacetime, \emph{Commun. Math. Phys.}
\textbf{50} (1984), 157--165.

\bibitem{brasil} A. Brasil and A.G. Colares, On constant mean curvature spacelike hypersurfaces in Lorentz manifolds, \emph{Mat. Contemp.}, \textbf{17} (1999) 99--136.



\bibitem{Ch} Y. Choquet-Bruhat,
Quelques propri\'et\'es des sousvari\'et\'es maximales d'une vari\'et\'e
lorentzienne, \emph{Cr. Acad. Sci. A Math. (Paris)  Serie A}, \textbf{281} (1975), 577--580


\bibitem{Ca} E. Calabi, Examples of Bernstein problems for some nonlinear
equations, \emph{P. Symp. Pure Math.}, \textbf{15} (1970),
223--230.

\bibitem{CST} E.J. Copeland, M. Sami and S. Tsujikawa,
Dynamics of dark energy, \emph{Int. J. Mod. Phys. D}, \textbf{15} (2006), 1753--1935.

\bibitem{CY} S.Y. Cheng and S.T. Yau, Maximal spacelike hypersurfaces in the
Lorentz-Minkowski spaces, \emph{Ann. of Math.},  \textbf{104} (1976),
407--419.

\bibitem{HE} S.W. Hawking and G.F.R. Ellis, \emph{The large scale structure of space-time}, Cambridge Univesity Press, London and New York, 1973. 

\bibitem{JKG} J.L. Jaramillo, J.A.V. Kroon and E. Gourgoulhon, From geometry to numerics: interdisciplinary aspects in mathematical and numerical relativity, \emph{Classical Quant. Grav.}, \textbf{25} (2008), 093001.

\bibitem{Lid} A.R. Liddle and D.H.  Lyth, \emph{Cosmological inflation and large-scale structure}, Cambridge University Press, 2000.

\bibitem{M-T} J.E. Marsden and F.J. Tipler, Maximal hypersurfaces and
foliations of constant mean curvature in General Relativity,
\emph{Phys. Rep.}, \textbf{66} (1980), 109--139.

\bibitem{N} S. Nishikawa, On maximal spacelike hypersurfaces in a Lorentzian manifold,
\emph{Nagoya Math. J.}, \textbf{95} (1984), 117--124.
.

\bibitem{O'N} B. O'Neill, \emph{Semi-Riemannian Geometry with
applications to Relativity}, Academic Press, 1983.



\bibitem{SP} S. Perlmutter et al.,  Measurements of Omega and Lambda from 42 High-Redshift Supernovae, \emph{Astrophysical J.}, \textbf{517} (1999), 565--586.



\bibitem{R} A.G. Riess, Observational evidence from supernovae for an accelerating universe and a cosmological constant, \emph{Astrophysical J.}, \textbf{116} (1999), 1009--1038.



\bibitem{RRS} A. Romero, R.M. Rubio and J.J. Salamanca,
Uniqueness of complete maximal hypersurfaces in spatially parabolic
generalized Robertson-Walker spacetimes, \emph{Classical Quant. Grav.}, \textbf{30} (2013) 115007--115020.


\bibitem{RRS2} A. Romero, R.M. Rubio and J.J. Salamanca, Complete maximal hypersurfaces in certain spatially open generalized Robertson–Walker spacetimes, \emph{Rev. R. Acad. Cienc. Exactas F\'is. Nat. Ser. A Math. RACSAM}, \textbf{109} (2015), 451--460.

\bibitem{Rub} R.M. Rubio, Complete constant mean curvature spacelike hypersurfaces in the Einstein-de Sitter spacetime, \emph{Rep. Math. Phys.}, \textbf{74} (2014), 127--133.

\bibitem{S} M. S\'anchez, On the geometry of generalized Robertson-Walker spacetimes: geodesics, \emph{Gen. Relat. Gravit.}, \textbf{30} (1998), 915--932.

\bibitem{Sha} G.S. Sharov, Multidimensional cosmological solutions of Friedmann type,
\emph{Theoret. and Math. Phys.}, \textbf{101} (1994), 1485--1491.






\end{thebibliography}
\end{document}